\documentclass[11pt, reqno]{amsart}
\usepackage{amssymb}

\keywords{Convex functions, approximately convex
functions, Hyers-Ulam Theorem, best constants}
\subjclass{Primary: 26B25, 41A44;  Secondary:  39B72,   51M16,  52A40}   
\usepackage{psfig}  

\theoremstyle{plain}

\newtheorem{theorem}{Theorem}

\newtheorem*{corollary}{Corollary}
\newtheorem{proposition}{Proposition}
\newtheorem{lemma}{Lemma}

\theoremstyle{remark}

\newtheorem*{remark}{Remark}
\newtheorem*{remarks}{Remarks}

\theoremstyle{definition}

\catcode`@=11 \@mparswitchfalse  
\newcounter{mnotecount}[section]

\renewcommand{\emptyset}{\varnothing}

\DeclareMathOperator{\sgn}{sgn}

\DeclareMathOperator{\supp}{supp}
\DeclareMathOperator{\Co}{Co}
\DeclareMathOperator{\sym}{Perm}
\newcommand{\cd}{,\dots,}
\newcommand{\R}{{\mathbb R}}		 
\newcommand{\N}{{\mathbb N}}	
\newcommand{\cn}{\colon}

\begin{document}
\numberwithin{equation}{section}   

\title[Extremal Approximately Convex Functions]
{Extremal Approximately Convex Functions and
the Best Constants in a Theorem of Hyers and Ulam}
\author{S. J. Dilworth}
\address{Department of Mathematics, University of South Carolina,
 Columbia, SC 29208, U.S.A.} \email{dilworth@math.sc.edu}

\author{Ralph Howard}
\address{Department of Mathematics, University of South Carolina,
 Columbia, SC 29208, U.S.A.} \email{howard@math.sc.edu}

\author{James W. Roberts}
\address{Department of Mathematics, University of South Carolina,
 Columbia, SC 29208, U.S.A.} \email{roberts@math.sc.edu}

\date{\today}

\begin{abstract} Let $n\ge1$ and $B\ge2$. A real-valued function $f$ 
defined on the $n$-simplex $\Delta_n$ is approximately convex with respect to 
$\Delta_{B-1}$ if
$$
f\left(\sum_{i=1}^B t_ix_i\right) \le \sum_{i=1}^B t_if(x_i) +1
$$
for all $x_1,\dots,x_B \in \Delta_n$ and all $(t_1,\dots,t_B)\in 
\Delta_{B-1}$.
We determine the extremal function of this type which vanishes on the 
vertices of $\Delta_n$.
We also prove a stability theorem of Hyers-Ulam type 
which yields as a special case 
the best constants in the Hyers-Ulam stability theorem for 
$\varepsilon$-convex functions.
 \end{abstract}


\maketitle

\section{Introduction} First we fix  some notation.  The standard $n$-simplex 
$\Delta_n$ is defined by 
\begin{equation*}
\Delta_n=
\Big\{(x(0),\dots,x(n)): \sum_{j=0}^n x(j)=1, x(j)\ge0, 0\le j\le n\Big\}. 
\end{equation*}
The vertices of $\Delta_n$ are denoted by $e(j)$ ($0\le j\le n$). For
$x \in \Delta_n$, the set $\{0\le j\le n\cn  x(j)\ne0\}$ is denoted by
$\supp x$.  Fix $B\ge 2$ and $n \ge 1$, and let $U$ be a convex subset
of $\mathbb{R}^n$. We say that a function $f\cn U\rightarrow\mathbb{R}$
is {\it approximately convex with respect to $\Delta_{B-1}$} if
\begin{equation*}
f\Big(\sum_{i=1}^B t_ix_i\Big) \le \sum_{i=1}^B t_if(x_i) +1 
\end{equation*}
for all $x_1,\dots,x_B \in U$ and all $(t_1,\dots,t_B)\in \Delta_{B-1}$.

In Section~2 we consider real-valued functions with domain $\Delta_n$
that are approximately convex with respect to $\Delta_{B-1}$. We show
that there exists an extremal such function satisfying the following:
(i) $E$ is approximately convex with respect to $\Delta_{B-1}$;
(ii) $E$ vanishes on the vertices of $\Delta_n$; (iii) if
$f\cn U\rightarrow\mathbb{R}$ is approximately convex with respect to
$\Delta_{B-1}$ and satisfies $f(e(j))\le 0$ for $j=0,\dots,n$, then
$f(x)\le E(x)$ for all $x\in\Delta_n$. Moreover, we obtain an explicit
formula for $E$, and we show that $E$ is concave and piecewise-linear
on $\Delta_n$ and continuous on the interior of $\Delta_n$. We also
calculate the maximum value of $E$.

In Section~3 we prove a stability theorem of Hyers-Ulam type for
approximately convex functions. In the case $B=2$, this result yields
the best constants in the well-known Hyers-Ulam stability theorem for
$\varepsilon$-convex functions \cite{HU}.

We refer the reader to the book \cite{HIR} for more information about
approximately convex functions and stability theorems.  Finally, for a
thorough treatment of extremal approximately midpoint-convex functions
and related results, we refer the reader to \cite{DHR1}.
\section{Extremal Approximately Convex Functions} 

Define a function $E\cn \Delta_n\to\mathbb{R}$ as follows (recall that
 $\sgn 0=0$ and $\sgn a=a/|a|$ if $a\ne0$):
\begin{equation} \label{eq: defofE}
E(x)=\min\Big\{\sum_{j=0}^n m(j)x(j): 
\sum_{j=0}^n \frac{\sgn x(j)}{B^{m(j)}} \le 
1,\  m(j)\ge0\Big\}.
\end{equation} 
If $x\in \Delta_n$ then $x(j)\ge 0$ and so $\sgn x(j)$ is either $0$
or $1$.  Note that if $A=\supp x$, then
\begin{equation} \label{eq: defofE2}
E(x)=\min\Big\{\sum_{j\in A} m(j)x(j): \sum_{j\in A} \frac{1}{B^{m(j)}}
\le 1,\  m(j)\ge0\Big\}.
\end{equation}

\begin{proposition} \label{prop: appconvex}
$E(e(j))=0$ for all $j$ and $E$ is approximately convex 
with respect to $\Delta_{B-1}$. 
\end{proposition}

\begin{proof} It is clear from \eqref{eq: defofE2} that $E(x) \ge 0$ for all 
$x$ and that $E(e(j))=0$ for all $j$.  Suppose that $x \in \Delta_n$
and that $x=\sum_{k=1}^Bt_kx_k$ for some $x_1,\dots,x_B\in
\Delta_n$. Let $A=\supp x$ and $A_k=\supp x_k$, and note that
$A\subseteq\bigcup_{k=1}^BA_k$.  For each $1\le k\le B$, we have
\begin{equation*} 
E(x_k)=\sum_{j\in A_k} m_k(j)x_k(j) 
\end{equation*}
for some $(m_k(j))_{j\in A_k}$ such that $\sum_{j\in A_k}1/B^{m_k(j)}
\le1$. For $j \in A$, let $C(j)=\{1\le k \le B: j \in A_k\}$
and let 
$$
M(j)= \min\{m_k(j): k \in C(j)\}.
$$
Note that 
\begin{equation*}
\frac{1}{B^{M(j)+1}}=\frac{1}{B}\frac{1}{B^{M(j)}} \le \frac{1}{B}
 \sum_{k \in C(j)}\frac{1}{B^{m_k(j)}}. 
\end{equation*}
Thus, 
\begin{equation*}\sum_{j\in A} \frac{1}{B^{M(j)+1}}\le
 \sum_{j\in A}\frac{1}{B}\sum_{k \in C(j)} 
\frac{1}{B^{m_k(j)}}\le\frac{1}{B}\sum_{k=1}^B\sum_{j\in 
A_k}\frac{1}{B^{m_k(j)}}\le1.
\end{equation*} 
Hence
\begin{align*} E\Big(\sum_{k=1}^Bt_kx_k\Big)=E(x)&\le\sum_{j\in A}(1+M(j))x(j)\\
&=\sum_{j\in A}(1+M(j))\sum_{k=1}^Bt_kx_k(j)\\
&= 1+ \sum_{k=1}^Bt_k\sum_{j\in A}M(j)x_k(j)\\
&= 1+ \sum_{k=1}^Bt_k\sum_{j\in A_k}M(j)x_k(j)\\
\intertext{(since $A_k \subseteq A$ if $t_k\ne0$)}
&\le 1+\sum_{k=1}^Bt_k\sum_{j\in A_k}m_k(j)x_k(j)\\
&=1+\sum_{k=1}^Bt_kE(x_k). 
\end{align*} 
Thus, $E$ is approximately convex with respect to $\Delta_{B-1}$. 
\end{proof}

\begin{lemma} \label{lem: lemma1}
If $m(j)\ge1$ for each $0\le j\le n$ and $\sum_{j=0}^n1/B^{m(j)}\le1$,
then $\{0,1,\dots,n\}$ is the disjoint union of sets $P_1,\dots,P_B$ such 
that 
\begin{equation*}
\sum_{j\in P_k}\frac{1}{B^{m(j)}}\le\frac{1}{B} 
\end{equation*}
for $k=1,\dots,B$.
\end{lemma} \begin{proof} Without loss of generality we may assume that
$1\le m(0)\le m(1)\le\dots\le m(n)$. We shall prove that the result holds for 
all $n\ge1$
by induction on $N=\sum_{j=0}^nm(j)$. Note that the result is vacuously true 
if $N=1$
and is trivial if $n\le B$. So suppose that $N\ge2$ and that $n>B$, so that 
$n-1> B-1\ge1$.
By inductive hypothesis, $\{0,1,\dots,n-1\}$ is the disjoint union of sets
$F_1,\dots,F_B$ such that
\begin{equation*}
\sum_{j\in F_k}\frac{1}{B^{m(j)}}\le\frac{1}{B}\end{equation*}
for $k=1,\dots,B$. Since $\sum_{j=0}^{n-1}1/B^{m(j)}<1$, and since $1\le 
m(0)\le m(1)\le\dots\le m(n)$, 
there exists $k_0$ such that
\begin{equation} \sum_{j\in F_{k_0}} \frac{1}{B^{m(j)}} \le 
\frac{1}{B}-\frac{1}{B^{m(n-1)}}
\le \frac{1}{B}-\frac{1}{B^{m(n)}}.
\end{equation} Put $P_{k_0}=P_{k_0}\cup\{n\}$ and $P_k=F_k$ for $k\ne k_0$ to 
complete the
induction. \end{proof}

\begin{theorem} \label{th: extremal}
$E$ is extremal, that is if $h\cn \Delta_n\to \R$ is approximately
convex with respect to $\Delta_{B-1}$ and $h(e(j))\le0$ for
$j=0,1,\dots,n$, then
$$h(x)\le E(x)\qquad \text{for all $x \in
\Delta_n$}.$$ 
\end{theorem}

\begin{proof} 
Let $s=|\supp x|$, so that $1\le s\le n+1$. The proof is by 
induction on $s$.
If $s=1$ then $x=e(j)$ for some $j$, so that  
$$E(x)=E(e(j))=0\ge h(e(j))=h(x).$$  As inductive hypothesis,
we suppose that $h(x)\le E(x)$ whenever $|\supp x|<s$. Now suppose that
$s\ge2$ and that $|\supp x|=s$.
Without loss of generality we may assume that $\supp x=\{0,\dots,s-1\}$,
so that
$E(x)=\sum_{j=0}^{s-1} m(j)x(j)$, where $\sum_{j=0}^{s-1} 1/{B^{m(j)}}\le1$.
Note that each $m(j) \ge1$ since $s \ge 2$.

If $\sum_{j=0}^{s-1} 1/{B^{m(j)}}\le1/B$, let $P_1=\{0,\dots,s-2\}$,
$P_2=\{s-1\}$, and $P_k=\emptyset$ for $2<k\le B$. Note that $|P_k|<s$
for $1\le k \le B$ and that $\sum_{j\in P_k} 1/B^{m(j)}\le 1/B$.

On the other hand, if
$\sum_{j=0}^{s-1} 1/{B^{m(j)}}>1/B$, then 
applying Lemma~\ref{lem: lemma1}
with $n=s-1$, we can write $\{0,1,\dots,s-1\}$ as the disjoint
union of sets $P_1, \dots,P_B$ such that $\sum_{j\in P_k} 1/B^{m(j)}\le 1/B$
for each $1\le k\le B$. Note that this implies that $|P_k|<s$ for $1 \le k \le B$.

If $P_k \ne \emptyset$, let 
$x_k = (1/t_k) \sum_{j\in P_k} x(j)e(j)$,
where $t_k=\sum_{j\in P_k}x(j)$.  If $P_k=\emptyset$, let $x_k=e(0)$ and let 
$t_k=0$. Thus $x=\sum_{k=1}^Bt_kx_k$, where $t_k\ge0$ and $\sum_{k=1}^B t_k=1$.
Note that $$|\supp x_k| = \max\{1, |P_k|\}<s \qquad(1\le k \le B).$$
If $P_k\ne\emptyset$, then $m(j)\ge 1$ for all $j\in P_k$, and $\sum_{j\in 
P_k}1/B^{m(j)-1}
\le 1$.  
Since
$|\supp x_k|<s$,  our inductive hypothesis implies that $h(x_k)\le 
E(x_k)$. Finally,
{\allowdisplaybreaks 
\begin{align*}
h(x) &= h\Big(\sum_{k=1}^Bt_kx_k\Big)
\le 1+\sum_{k=1}^Bt_kh(x_k)
\le1+\sum_{P_k\ne\emptyset} t_kE(x_k)\\
&\le1+ \sum_{P_k\ne\emptyset} t_k\sum_{j\in P_k}(m(j)-1)x_k(j)\\
&=1+\sum_{P_k\ne\emptyset} \sum_{j\in P_k}(m(j)-1)x(j)\\
&=1+\sum_{j=0}^{s-1}m(j)x(j)-\sum_{j=0}^{s-1}x(j)\\
&=\sum_{j=0}^{s-1}m(j)x(j) 
=E(x). 
\end{align*} 
} 
This completes the induction.
\end{proof}
Following the convention that $x\log_Bx=0$ when $x=0$,
the \textit{entropy} function $F\colon\Delta_n \rightarrow \mathbb{R}$
is defined as follows: 
\begin{equation*}
F(x) = -\sum x(j) \log_B x(j).\end{equation*}
\begin{proposition} $F$ is approximately convex with respect to $\Delta_{B-1}$
and satisfies
$$F(x) \le E(x) \le F(x)+1\qquad(x \in \Delta_n).$$ 
\end{proposition}
\begin{proof} Let $x \in \Delta_n$. 
A standard Lagrange multiplier calculation yields
\begin{equation}
 \label{eq: defofF}
F(x)=\min\Big\{\sum_{j\in A} y(j)x(j): \sum_{j\in A} \frac{1}{B^{y(j)}}
\le 1,\  y(j)\ge0\Big\},
\end{equation} where $A=\supp x$.
Using \eqref{eq: defofF} in place of  \eqref{eq: defofE2}, minor changes in
the proof of Proposition~\ref{prop: appconvex} 
show that $F$ is approximately convex with respect to $\Delta_{B-1}$. 
Suppose that \begin{equation}
F(x)= \sum_{j\in A} y(j)x(j) \end{equation}
for some $y(j)\ge0$ satisfying $\sum_{j\in A} 1/B^{y(j)} \le 1$. Let
$m(j)=\lceil y(j) \rceil$. Then $\sum_{j\in A} 1/B^{m(j)} \le 1$, and so
\begin{equation*}
E(x) \le \sum_{j\in A} m(j)x(j) \le \sum_{j \in A} (y(j)+1)x(j) = F(x)+1.
\end{equation*} On the other hand, since $F$ is approximately convex
 with respect to $\Delta_{B-1}$, it follows from Theorem~\ref{th: extremal} that
$F(x) \le E(x)$. \end{proof}

\begin{proposition}\label{E-basics} 
(i) $E$ is piecewise-linear and  the restriction of $E$
to each open
 facet of $\Delta_n$ is continuous. \newline
(ii) $E$ is lower semi-continuous; \newline 
(iii) $E$ is concave.
\end{proposition}

\begin{figure}[ht]
\centering
\mbox{\psfig{file=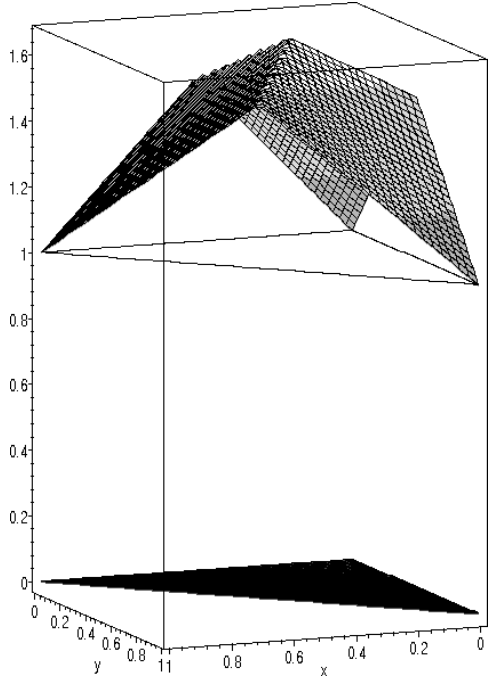,width=2.5in}}
\caption[]{\footnotesize Graph of $y=E(x,y,1-x-y)$ for $B=1$
over the simplex $0\le y\le 1-x\le1$ showing the discontinuity along the
boundary. On the boundary $E_S$ has the value~$1$ except at the three
vertices where it has the value~$0$.}
\label{Esimp-graph}
\end{figure}

\begin{proof} 
To prove that $E$ is piecewise linear it is enough to show that $E$
is piecewise linear on the interior $\Delta_n^\circ$ of $\Delta_n$.
For then by an induction on $n$ we will have that $E$ is piecewise
linear on $\Delta_n^\circ$ and the induction hypothesis implies that it is
piecewise linear when restricted to any of the facets of $\Delta_n$,
which implies that $E$ is piecewise linear on $\Delta_n$.
For fixed $n$ and $B$ let
$$
\mathcal{F}(n,B):=\left\{(m_0\cd m_n) : m_k\in \N,\  
\sum_{k=0}^n\frac{1}{B^{m_k}}\le 1\right\}
$$
be the set of feasible $(n+1)$-tuples.  For $(m_0\cd m_n)\in
\mathcal{F}(n,B)$ let $\Lambda_{(m_0\cd m_n)}\Delta_n\to \R$ be the
linear function
$$
\Lambda_{(m_0\cd m_n)}(x_0\cd x_n)=m_0x_0+m_1x_1+\cdots+m_nx_n
$$
so that $E\cn \Delta_n\to \R$ is given by 
$$
E(x)=\min\{\Lambda_{(m_0\cd m_n)}(x):(m_0\cd m_n)\in
\mathcal{F}(n,B)\}.
$$
Let
\begin{align*}
\mathcal{E}(n,B) :=\{(m_0\cd m_n)&\in \mathcal{F}(n,B):\\
 &\Lambda_{(m_0\cd m_n)}(x)=E(x) \text{ for some } x\in \Delta_n^\circ\}
\end{align*}
be the set of extreme $(n+1)$-tuples.  Then 
$$
E\big|_{\Delta_n^\circ}(x)=
\min\{\Lambda_{(m_0\cd m_n)}(x):(m_0\cd m_n)\in \mathcal{E}(n,B)\}
$$
and therefore showing that $E\big|_{\Delta_n^\circ}$ is piecewise
linear is equivalent to showing that $\mathcal{E}(n,B)$ is finite.

\begin{lemma}\label{E-min}
Let $(m_0\cd m_n)\in \mathcal{E}(n,B)$ and $(m_0'\cd m_n')\in
\mathcal{F}(n,B)$ with $m_k'\le m_k$ for $0\le k\le n$.  Then $(m_0'\cd
m_n')=(m_0\cd m_n)$.
\end{lemma}

\begin{proof}
For if not then there is an index $k$ with $m_k'<m_k$.  As all the
components of $x=(x_0\cd x_n)$ are positive on $\Delta_n^\circ$ this
implies that on $x\in \Delta_n^\circ$
\begin{align*}
E(x)&\le \Lambda_{(m_0'\cd m_n')}(x)=\Lambda_{(m_0\cd m_n)}(x)+
\Lambda_{(m_0'\cd m_n')}(x)-\Lambda_{(m_0\cd m_n)}(x)\\
&\le  \Lambda_{(m_0\cd m_n)}(x) +(m_k'-m_k)x_k<\Lambda_{(m_0\cd m_n)}(x).
\end{align*}
This contradicts that for $(m_0\cd m_n)\in\mathcal{E}(n,B)$ there is
an $x\in \Delta_n^\circ$ with $\Lambda_{(m_0\cd m_n)}(x)=E(x)$.
\end{proof}

Let $\sym(n+1)$ be the group of permutations of $\{0,1\cd n\}$.  Then
it is easily checked that $\mathcal{E}(n,B)$ is invariant under the
action of $\sym(n+1)$ given by $\sigma(m_0,m_1\cd
m_n)=(m_{\sigma(0)},m_{\sigma(1)}\cd m_{\sigma (n)})$.  Therefore if
$\mathcal{E}^*(n,B)$ is the set of monotone decreasing elements of
$\mathcal{E}(n,B)$, that is
$$
\mathcal{E}^*(n,B):=\{(m_0\cd m_n)\in \mathcal{E}(n,B): m_0\ge m_1\ge
\cdots \ge m_n\},
$$
then 
$$
\mathcal{E}(n,B)=\{\sigma(m_0\cd m_n): (m_0\cd m_n)\in
\mathcal{E}^*(n,B), \sigma\in \sym(n+1)\}
$$
and to show that $\mathcal{E}(n,B)$ is finite it is enough to show
that $\mathcal{E}^*(n,B)$ is finite.  

\begin{lemma}\label{ex-C}
Suppose that $n\ge0$. 
Let $m_0\ge m_1\ge \cdots \ge m_n$ be a non-increasing sequence of $(n+1)$
positive integers, and let $C$ be a positive real number such that 
$$
\sum_{k=0}^n\frac{1}{B^{m_k}}\le C,
$$
and such that if $m_0',m_1'\cd m_n'$ are any positive integers with
$m_k'\le m_k$ for $0\le k\le n$, then 
$$
\sum_{k=0}^n\frac{1}{B^{m_k'}} \le C
$$
implies that $(m_0'\cd m_n')=(m_0\cd m_n)$.  
(We will say that \emph{$(m_0\cd m_n)$
is  extreme for $(n,C)$}.)  Let
$$
\eta=\eta(n,C):=\min\{j\ge 2: CB^j\ge n+B\}.
$$
Then $m_n<\eta(n,C)$.  (The explicit value of $\eta$ is
$\eta(n,C)=\max\{2,\lceil\log_B((n+B)/C)\rceil\}$.)
\end{lemma}

\begin{proof}
 From the definition of $\eta$ we have $\eta\ge2$ and $CB^\eta\ge n+B$
which is equivalent to
$$
\frac{n+1}{B^\eta}\le C-\frac{1}{B^{\eta-1}}+\frac{1}{B^\eta}.
$$
Assume, toward a contradiction, that $m_n\ge \eta$.  Then
$$
\frac{1}{B^{m_0}}+\cdots+\frac{1}{B^{m_{n-1}}}+\frac{1}{B^{m_n}}
\le \frac{n+1}{B^\eta}\le C-\frac{1}{B^{\eta-1}}+\frac{1}{B^\eta}.
$$
This can be rearranged to give
$$
\frac{1}{B^{m_0}}+\cdots+\frac{1}{B^{m_{n-1}}}+\frac{1}{B^{\eta-1}}\le
C+\frac{1}{B^\eta}-\frac{1}{B^{m_n}}\le C.
$$
This contradicts that $(m_0\cd m_n)$ is $(n,C)$ extreme and completes
the proof.
\end{proof} 
We now prove $\mathcal{E}^*(n,B)$ is finite.  First some notation.
For positive integers $l_1\cd l_j$ let $C(l_1\cd
l_j):=1-\sum_{i=1}^j1/B^{l_j}$.  If $(m_0\cd m_n)\in
\mathcal{E}^*(n,B)$ then by Lemma~\ref{E-min} (and with the
terminology of Lemma~\ref{ex-C}) for each $j$ with $1\le j\le n$ the
tuple $(m_0\cd m_{n-j})$ is $(n-j,C(m_{n-j+1}\cd m_n))$ extreme, and
$(m_0\cd m_n)$ itself is $(n,1)$ extreme.  Therefore, by
Lemma~\ref{ex-C}, $m_n < \eta(n,1)$, whence there are only a finite
number of possible choices for $m_n$.  For each of these choices of
$m_n$ we can use Lemma~\ref{ex-C} again to get $m_{n-1}<
\eta(n-1,C(m_n))$, and so there are only finitely many choices for the
ordered pair $(m_{n-1},m_n)$.  And for each of these pairs
$(m_{n-1},m_n)$ we have that so there are only finitely many
possibilities for $m_{n-2}$.  Continuing in this manner it follows
that $\mathcal{E}^*(n,B)$ is finite.  This completes the proof that
$E_S^{\Delta_n}$ is piecewise linear and thus point~(i) of
Propsition~\ref{E-basics}

To prove
point~(ii) let  $A$ be a nonempty subset of $\{0,1,\dots,n\}$. In
proving point~(i) we have seen that there is a finite collection
$\mathcal{L}(A)$ of linear mappings $\Lambda\cn
\Delta_n\rightarrow\mathbb{R}$, each one of the form
$\Lambda(x)=\sum_{j\in A} m(j)x(j)$ for some nonnegative integers
$m(j)$, $j=0,1,\dots,n$, with $\sum_{j\in A}1/B^{m(j)}\le1$, such that
\begin{equation}\label{eq: defofE3}
E(x)= \min\{\Lambda(x): \Lambda\in \mathcal{L}(A)\} 
\end{equation}
for all $x \in \Delta_n$ such that $\supp x = A$.  Clearly, we may
also assume that $\mathcal{L}(B)\subseteq \mathcal{L}(A)$ whenever $A
\subseteq B$.  Suppose that $(x_i)_{i=1}^\infty \subseteq \Delta_n$
and that $x_i \rightarrow x$ as $i\rightarrow\infty$. Note that
$\supp x \subseteq \supp x_i$ for all sufficiently large $i$, so that
$\mathcal{L}(\supp x_i)\subseteq \mathcal{L}(\supp x)$ for all
sufficiently large $i$.  Thus,
\begin{align*} E(x) &= \min\{T(x): T\in \mathcal{L}(\supp x)\}\\
&= \lim_{i\rightarrow\infty}\min\{T(x_i):T\in\mathcal{L}(\supp x)\}\\
&\le \liminf_{i\rightarrow\infty}\min\{T(x_i):T\in\mathcal{L}(\supp x_i)\}\\
&=\liminf_{i\rightarrow\infty}E(x_i). 
\end{align*}
Thus, $E$ is lower semi-continuous. 

Finally we prove point~(iii).  It follows from \eqref{eq: defofE3}
that the restriction of $E$ to the interior of any facet is the
minimum of a finite collection of linear functions, and hence is
continuous and concave. The lower semi-continuity of $E$ forces $E$ to
be concave on all of $\Delta_n$. 
\end{proof}

\begin{remark}\label{E*values}
The algorithm implicit in the proof that $\mathcal{E}^*(n,B)$ is finite
is rather effective for small values of $n$.  In the case of most
interest, when $B=2$ so that $S=\Delta_1$, it can be used to show
\begin{align*}
\mathcal{E}^*(2,2)&=\{(2,2,1)\},\qquad 
\mathcal{E}^*(3,2)=\{(3,3,2,1),(2,2,2,2)\}\\
\mathcal{E}^*(4,2)&=\{(4,4,3,2,1), (3,3,2,2,2)\},\\
\mathcal{E}^*(5,2)&=\{5,5,4,3,2,1), (3,3,3,3,2,2)\}.
\end{align*}
When $n=2$ this leads to the explicit formula
$$
E(x,y,1-x-y)=\min\{1+x+y, 2-x, 2-y\}
$$
for $0<x<1-y<1$. (Cf.~Figure~\ref{Esimp-graph}).  The sets
$\mathcal{E}^*(n,2)$ can be used to give messier, but equally explicit
formulas, for higher values of $n$.\qed
\end{remark}

\begin{proposition} \label{prop: maximum}
The maximum of $E$ is given by 
\begin{equation} \label{eq: kappa(B,n)}
\kappa(n,B) = \lfloor \log_Bn\rfloor +\frac{\lceil B(n+1-B^{\lfloor 
\log_Bn\rfloor})/(B-1)\rceil}
{n+1} 
\end{equation}
\end{proposition} 

For small values of $B$ and $n$, $\kappa_S(n)$ is given in
Table~\ref{kappa-table}.
\begin{table}[hb]
{\tiny
$$
\begin{array}{c|ccccccccccc}
B\backslash n& 1&2&3&4&5&6&7&8&9&10\\ \hline
2&1.0& 1.6667 & 2.0000 & 2.4000 & 2.6667 & 2.8571 & 3.0000 & 3.1111 & 
3.4000 & 3.5455 \\
3&1.0& 1.0& 1.5000 & 1.6000 & 1.8333 & 1.8571 & 2.0000 & 2.0000 & 2.2000 & 
2.2727 \\
4&1.0& 1.0& 1.0& 1.4000 & 1.5000 & 1.5714 & 1.7500 & 1.7778 & 
1.8000 & 1.9091 \\
5&1.0& 1.0& 1.0& 1.0& 1.3333 & 1.4286 & 1.5000 & 1.5556 & 1.7000 & 
1.7273 \\
6&1.0& 1.0& 1.0& 1.0& 1.0& 1.2857 & 1.3750 & 1.4444 & 1.5000 & 
1.5455 \\
7&1.0& 1.0& 1.0& 1.0& 1.0& 1.0& 1.2500 & 1.3333 & 1.4000 & 1.4545
 \\
8&1.0& 1.0& 1.0& 1.0& 1.0& 1.0& 1.0& 1.2222 & 1.3000 & 1.3636 \\
9&1.0& 1.0& 1.0& 1.0& 1.0& 1.0& 1.0& 1.0& 1.2000 & 1.2727 \\
10&1.0& 1.0& 1.0& 1.0& 1.0& 1.0& 1.0& 1.0& 1.0& 1.1818 \\
11&1.0& 1.0& 1.0& 1.0& 1.0& 1.0& 1.0& 1.0& 1.0& 1.0
\end{array}
$$}
\caption[]{\footnotesize
Values of $\kappa(n,B)$ for $2\le B\le 11$
and $1\le n\le 10$.}
\label{kappa-table}
\end{table}

\begin{proof}
$E$ is a symmetric function of $x(0),\dots,x(n)$ and $E$ is also concave.
 Thus $E$ achieves its maximum at the barycenter $\overline{x}=
(1/(n+1))\sum_{j=0}^ne(j)$.  So there exist nonnegative
 integers $m(j)$ ($j=0,1,\dots,n$) such that
$E(\overline{x})=(1/(n+1))\sum_{j=0}^n m(j)$ and $\sum_{j=0}^n 1/B^{m(j)} \le 
1$.
We may also assume that $(m(j))_{j=0}^n$ have been chosen to minimize 
$\sum_{j=0}^n 1/B^{m(j)}$
among all  possible choices of $(m(j))_{j=0}^n$. 
Suppose that there exist $i$ and $k$ such that $m(k)\ge m(i)+2$. Note that
\begin{equation} \frac{1}{B^{m(i)+1}}+\frac{1}{B^{m(k)-1}}
\le\frac{2}{B^{m(i)+1}}\le 
\frac{B}{B^{m(i)+1}}<\frac{1}{B^{m(i)}}+\frac{1}{B^{m(k)}}.
\end{equation} Thus replacing $m(i)$ by $m(i)+1$ and replacing $m(k)$ by 
$m(k)-1$
leaves $(1/(n+1))\sum_{j=0}^n m(j)$ unchanged while it reduces $\sum_{j=0}^n 
1/B^{m(j)}$, which
contradicts the choice of  $(m(j))_{j=0}^n$. Thus $|m(i)-m(k)|\le1$
for all $i,k$. It follows that there exist integers $\ell\ge0$ and $1\le s\le 
n+1$
such that 
\begin{equation} \label{eq: kappa}
\kappa(n,B)= \frac{\ell(n+1-s)+(\ell+1)s}{n+1}=\ell+\frac{s}{n+1} 
\end{equation}
and 
\begin{equation} \label{eq: kappa2}
 \frac{n+1-s}{B^\ell}+\frac{s}{B^{\ell+1}} \le 1. \end{equation}
 Moreover, it is clear from \eqref{eq: kappa} that
$\ell$ is the least nonnegative integer satsifying \eqref{eq: kappa2}
for some $1\le s \le n+1$, i.e. \begin{equation*}
\ell=\lfloor \log_Bn\rfloor. \end{equation*}
For this value of $\ell$ it is 
clear from \eqref{eq: kappa} that $s$ is the smallest integer in the range $1 
\le s \le
n+1$ satisfying \eqref{eq: kappa2}, i.e.
\begin{equation*}
s=\left\lceil\frac{B(n+1)-B^{\ell+1}}{B-1}\right\rceil
=\left\lceil\frac{B}{B-1}(n+1-B^\ell)\right\rceil. 
\end{equation*}
Substituting these values for $\ell$ and $s$ into \eqref{eq: kappa} gives 
\eqref{eq: kappa(B,n)}.
\end{proof}

\section{Best Constants in Stabilty Theorems of Hyers-Ulam Type}
Hyers and Ulam \cite{HU} introduced the following definition. Fix 
$\varepsilon>0$.
A function $f\cn U\rightarrow\mathbb{R}$, where $U$ is a convex subset of 
$\mathbb{R}^n$,
 is $\varepsilon$-convex if \begin{equation*}
 f(tx+(1-t)y)\le tf(x)+(1-t)f(y)+\varepsilon \end{equation*}
for all $x,y\in U$ and all $t\in[0,1]$.

Note that $f$ is $\varepsilon$-convex 
if and only if
$(1/\varepsilon)f$ is approximately convex with respect to $\Delta_1$. 
So let us generalize this notion by defining 
$f$ to be \textit{$\varepsilon$-convex with respect to
$\Delta_{B-1}$} if $(1/\varepsilon)f$ is approximately 
convex with respect to $\Delta_{B-1}$.

The proof of the following theorem is adapted from Cholewa's
proof \cite{C} of the Hyers-Ulam stability theorem for $\varepsilon$-convex 
functions.

\begin{theorem} \label{th: HyersUlam1}
 Suppose that $U\subseteq \mathbb{R}^n$ is convex and that 
$f\cn U\rightarrow\mathbb{R}$ is $\varepsilon$-convex with respect to 
$\Delta_{B-1}$.
Then there exist convex functions $g,g_0\cn U\rightarrow\mathbb{R}$ such that
\begin{equation*}
g(x)\le f(x)\le g(x)+ \kappa(n,B)\varepsilon\qquad\text{and}\qquad 
|f(x)-g_0(x)|\le\frac{\kappa(n,B)\varepsilon}{2}
\end{equation*} for all $x\in U$. Moreover, $\kappa(n,B)$ is the best constant
in these inequalities. 
\end{theorem}

\begin{proof}
By replacing $f$ by $f/\varepsilon$, we may assume that
$\varepsilon=1$.  Set $W=\{(x,y)\in U\times\mathbb{R}: y\ge
f(x)\}\subseteq\mathbb{R}^{n+1}$ and define $g$ by
\begin{equation}
g(x)=\inf\{y:(x,y)\in \Co(W)\}. 
\end{equation}
Clearly $-\infty \le g(x) \le f(x)$.
Suppose that $(x,y) \in \Co(W)$. By Caratheodory's Theorem
(see e.g. \cite[Thm.\ 17.1]{R}) there exist $n+2$
points $(x_0,y_0),\dots,(x_{n+1},y_{n+1})\in W$ such that $(x,y) \in \Delta 
:= 
\Co(\{(x_0,y_0),\dots,(x_{n+1},y_{n+1})\})$. Let 
$\overline{y}=\min\{\eta:(x,\eta)\in\Delta\}$.
Then $(x,\overline{y})$ lies on the boundary of $\Delta$ and so it is a 
convex combination of $n+1$ of the points $(x_0,y_0),\dots,(x_n,y_n)$.
Without loss of generality,
$(x,\overline{y})
= \sum_{j=0}^n t_j(x_j,y_j)$ for some $(t_0,\dots,t_n)\in \Delta_n$.
Note that 
\begin{equation*}
h\Big(\sum_{j=0}^n x(j)e(j)\Big):= f\Big(\sum_{j=0}^n x(j) 
x_j\Big)-
\sum_{j=0}^n x(j) f(x_j)\qquad(x \in \Delta_n)
\end{equation*}
 is approximately convex
with respect to $\Delta_{B-1}$ and satisfies $h(e(j))=0$ for $j=0,1,\dots, n$.
By Proposition~\ref{prop: maximum}, $\max_{x\in\Delta_n} h(x) \le 
\kappa(n,B)$. Thus
\begin{align*}
 y &\ge \overline{y} = \sum_{j=0}^n t_j y_j
 = \sum_{j=0}^n t_j f(x_j)\\
&= f\Big(\sum_{j=0}^n t_j x_j\Big)-h\Big(\sum_{j=0}^n t_je(j)\Big)\\
&\ge f\Big(\sum t_j x_j)\Big)- \kappa(n,B)\\ &= f(x)- \kappa(n,B). 
\end{align*}
Taking the infimum over all $y$ yields $g(x) \ge f(x)- \kappa(n,B)$,
 i.e. $f(x) \le g(x)+ \kappa(n,B)$. Finally, set $g_0(x) = g(x)+ 
\kappa(n,B)/2$.

The fact that $\kappa(n,B)$ is the best constant follows by taking
$f$ to be $E$, where $E$ is the extremal approximately convex function 
(with respect to $\Delta_{B-1}$) with domain $\Delta_n$.
\end{proof}

Thus, setting $B=2$ in Theorem~\ref{th: HyersUlam1} gives the best
constants in the Hyers-Ulam stability theorem for $\varepsilon$-convex
functions \cite{HU}.

\begin{corollary} \label{eps-thm} 
Suppose that $U\subseteq \mathbb{R}^n$ is convex and that
$f\cn U\rightarrow\mathbb{R}$ is $\varepsilon$-convex.
Then there exist convex functions $g,g_0\cn U\rightarrow\mathbb{R}$ such that
\begin{equation*}
g(x)\le f(x)\le g(x)+ \kappa(n)\varepsilon\qquad\text{and}\qquad 
|f(x)-g_0(x)|\le\frac{\kappa(n)\varepsilon}{2}
\end{equation*} for all $x\in U$, where \begin{equation*}
\kappa(n)= \lfloor \log_2n\rfloor+\frac{2(n+1-2^{\lfloor 
\log_2n\rfloor})}{n+1}.
\end{equation*}
Moreover, $\kappa(n)$ is the best constant
in these inequalities. 
\end{corollary}

\begin{remarks}
 1. The value $\kappa(2)=5/3$ was first obtained by Green 
\cite{G}.
 The value $\kappa(2^n-1)=n$ was obtained by a different argument in 
\cite{DHR2}.

2. Note that $\kappa(3)=2$, $\kappa(4)=12/5$, $\kappa(5)=8/3$,
$\kappa(6)=20/7$, $\kappa(7)=3$, etc. These values improve the
constants obtained by Cholewa~\cite{C}.

3. The best constants corresponding to $\kappa(n)$ 
 for approximately midpoint-convex functions 
 were obtained
in~\cite{DHR1}.
\end{remarks}

\end{document}